# The mathematics in middle-aged Arab caliphate and it application to contemporary teaching in high schools

## Adnan Sharaf Ali[1], Stancho Pavlov[2], Krasimir Yordzhev[3]

[1]engadnansharaf@yahoo.com, [2]stancho_pavlov@yahoo.com, [3]yordzhev@swu.bg

[1,3]South-West University "Neofit Rilski", Blagoevgrad, Bulgaria

[2]Asen Zlatarov University, Burgas, Bulgaria

**Abstract:** This paper offers a glimpse of the major contributions made by Arabs to mathematics in middle ages history period. Its purpose is to stimulate interest in an object based on mutual respect and understanding. We give a short list of the most popular Arabic mathematicians and we comment some of their studies for using them in the high school teaching.

### 1. Reasons for the flourishing of science in Arabic State in middle ages

**1.1. Common language.** Arabic is the language of more than 240 million inhabitants of the Arabic world. Arabic language emerged from ancient Semitic traditions that were steeped in poetry over thousands of years before it became Arabic. The first mention of Arabic language is in the Qur'an. Qur'anic Arabic was the first Semitic language that codified its phonetic presentation (a complete recording of its sound and rhythm), and brought linguistic stability to Semitic. As the Empire spread, the Arabic language and culture was enriched by contacts with other civilizations. A great translation activities was in force, in which many ancient tracts where translated into Arabic. Then there was a proverb "It is better to translate one good book than to write a hundreds bad".

**1.2. Relationship with other cultures.** As the empire spread, the Arabic language and culture was enriched by contacts with other civilizations – Greeks, Indians, Persians, Romans Copts and Chinese. When in 1492 the Christians entered the libraries of Cordova they was fascinated by the thousands books in them. At his time Christians, libraries possessed about ten books. The European scientific world must be thankful to Arabian scholars for preserving writings of Euclid, Archimedes, and Aristotle and Claudius Ptolemy. They are translated from Arabic to Latin. The Arab scientist Muhammad ibn Musa al-Khwarizmi introduced the Hindu-Arabic numeral system in Arabic world and later it spreads over all Europe. [5]

**1.3. Common religion.** Islam originated in the Arabian Peninsula in 622 A.D. The Five Pillars of Islam express the Islamic acts:

(1) *Shahada* requires confess the unity of God and Muhammad. This involves the repetition of the formula "There is no God but Allah and Muhammad is the messenger of Allah."

(2) *Sala* Prayer is required five times a day.

(3) *Zaka* embodies the principle of social responsibility. What belongs to the believer also belongs to the community. By donation a proportion of his wealth for public use does believer legitimize what he retains.

(4) *Saum* means observing Ramadan (fasting month), the month during which God send the Qur'an to Gabriel and he reviled it to Muhammad.

(5) *Al-Hajj Pilgrimage* to Mecca. All pilgrimage from various classes express the full equality seeking to gain the favor of God.

The Islamic community throughout the world id united by two believes:

(1) Oneness of the God
(2) The mission of his Prophet [1]

**1.4. Science and education.** Among the early elementary educational institutions where the mosque schools. Arabs used basic teachings as a starting point from which to begin a revolution in educational during Abbasids dynasty (750 -1258 A.D.). In X century in Baghdad were 3 000 mosques. In XIV century Alexandria had 12 000 mosques, all of which played an important role in education. Students could attend several classes a day, sometimes traveling from one mosque to another. There was a strict discipline in these schools. The students respected teachers and there were formal rules of behavior. Laughing, talking, joking or disrespectful behavior of any kind were not permitted. The madrasa gave rise to many universities in the Arab caliphate Al-Azhar University in Cairo preceded other universities in Europe by two centuries. Traveling to other cities to seek knowledge under the direction of different masters was a common practice in the Arab Caliphate. From Kurasan to Egypt, to West Africa and Spain students and teachers journeyed to attend different classes.

**1.5. Military forces.** The Arab army maintained a high level of discipline, strategic prowess and organization. In its time, the army was a powerful and very effective force.
Its size was initially 13,000 troops in 632, but as the Caliphate expanded, the army gradually grew to 100,000 troops by 657. Only Muslims were allowed to join the army as regular troops. However, during the Muslim conquest of Roman Syria (633-638) some 4,000 Byzantine soldiers under their commander Joachim converted to Islam and served as troops. [6]

**2. Islamic Mathematics and Contribution to the civilization**

**2.1** We will begin with three old puzzles. [2]

The following problem was invented by an unknown arabic mathematician many thousand years ago:

A rich old Arab has three sons. When he died, he willed his 17 camels to the sons, to be divided as follows: First Son to get 1/2 of the camels. Second Son to get 1/3 of the camels. Third son to get 1/9 of the camels. The sons are sitting there trying to figure out how this can possibly be done, when a very old wise man goes riding by on his camel. They stop him and ask him to help them solve their problem. Without hesitation, he divides the camels properly and continues riding on his way. How did he do it? The old man temporarily added his camel to the 17, making sum of camels equal to 18.
$$18/2 = 9; 18/3 = 6; 18/9 = 2$$
$$9+6+2=17$$
He then takes his camel back and rides away.

The second one:

Hunter met two shepherds, one of whom has three slices of bread, and the another has five. All the pieces are of equal size. The three men shared the whole bread. After eating, the hunter gave shepherds eight coins. How the shepherds have to share the money?

Answer 1:7.

The third one: (The **Al-Kashi**'s task)

Pay for an employee work for a month (that is thirty days) is 10 dinars and a dress. The employee worked three days only and earned the dress only. What is the price of the dress?

Answer 10/9.

## 2.2. Abu Abdallah Muhammad ibn Musa al-Khwarizmi (790 – 850) Baghdad

In his book, "The Compendious Book on Calculation by Completion and Balancing" Al-Khwarizmi deals with ways to solve for the positive roots of the second and third degree polynomial equations. He also in his algebra was rhetoric, which means that the equations were written out in full sentences. The translations to symbolic algebra, where only symbols are used, can be seen in the work of Ibn al-Bunna' al Marakushi.

Al-Khwarizmi introduced the decimal place value system and make it widespread throughout Arabic and then in Europe. Al-Khwarizmi demonstrates also the basic operations of addition, subtraction, division and multiplication with decimal numbers. [10]

The word *algebra* comes directly from his book**, "Hisab al-jabr** w'al-muqabala" (**"The Science of Restoration and Reduction"**)**. The algorithm are derived from the name of al-Khwarizmi.

## 2.3. Abu'l-Hasan Thabit ibn Qurra (826 - 901) Baghdad

He made many contributions to mathematics, particularly in number theory for amicable and perfect numbers. He discovered a theorem which allowed pairs of amicable numbers to be found, that is two numbers such that each is the sum of the proper divisors (greater than zero and less than the number itself) of the other. The theorem states that if $p = 3.2^{n-1} - 1$, $q = 3.2^n - 1$ and $r = 9.2^{2n-1} - 1$ are primes then $2^n.p.q$ and $2^2.r$ are amicable.

Euler generalizes this theorem in 1750.

**Al-Baghdadi** (born 980) looked at a slight variant of Thabit ibn Qurra's theorem, while **al-Haytham** (born 965) seems to have been the first to attempt to classify all even perfect numbers (numbers equal to the sum of their proper divisors) as those of the form $2^{k-1}(2^k - 1)$ where $2^k - 1$ is prime.

Continuing the story of amicable numbers, from which we have taken a diversion, it is worth noting that they play a large role in Arabic mathematics. **Al-Farisi** (born 1260) gave a new proof of Thabit ibn Qurra's theorem, introducing important new ideas concerning factorisation and combinatorial methods. He also gave the pair of amicable numbers 17296 and 18416 which have been attributed to Euler, but we know that these were known earlier than al-Farisi, perhaps even by Thabit ibn Qurra himself. Although outside our time range for Arabic mathematics in this article, it is worth noting that in the 17[th] century the Arabic mathematician **Mohammed Baqir Yazdi** gave the pair of amicable number 9,363,584 and 9,437,056 still many years before Euler's contribution.

## 2.4. Thabit ibn Qurra generalized Pythagoras's theorem to an arbitrary triangle, as did Pappus.

**Theorem** (Generalization of Pythagorean Theorem) Construct two points A$^/$ and B$^/$ on the side BC of the triangle ABC so that ∠AA$^/$C = ∠BB$^/$C = γ. Then a$^2$+ b$^2$ = c(AA$^/$ + AA$^/$ ).

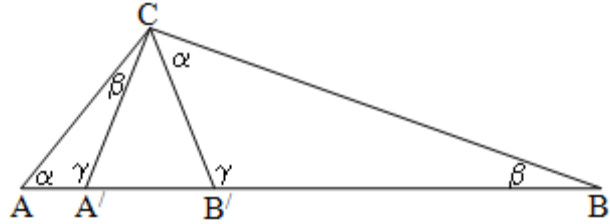

From the similarity of the triangles ABC and CB'B we deduce that AB:BC=BC:BB'.

$$\triangle ABC \approx \triangle CB'B \implies \frac{AB}{BC}=\frac{BC}{BB'} \implies cBB'=a^2$$

Analogically we obtain $cAA'=b^2$.
From these two equations we derive the claimed by adding.

### 2.5. Abu Ali al-Hasan ibn al-Hasan ibn al-Haytham (965 – 1039) Basra, he lived mainly in Cairo, Egypt

Al-Haytham, is the first mathematician, that we know to state Wilson's theorem, namely that if *p* is prime then 1+(*p*-1)! is divisible by *p*. It is unclear whether he knew how to prove this result. It is called *Wilson's theorem* because of a comment made by Waring in 1770 that John Wilson had noticed the result. There is no evidence that John Wilson knew how to prove it and most certainly Waring did not. Lagrange gave the first proof in 1771 and it should be noticed that it is more than 750 years after al-Haytham before number theory surpasses this achievement of Arabic mathematics.

For his volume computations, al-Haytham needed formulas for the sums of the first n integer cubes and the first n fourth powers. He may have used a diagram like that to describe the relationship

$$(4+1)\sum_{i=1}^{4}i = \sum_{i=1}^{4}i^2 + \sum_{p=1}^{4}\left(\sum_{i=1}^{p}i\right)$$

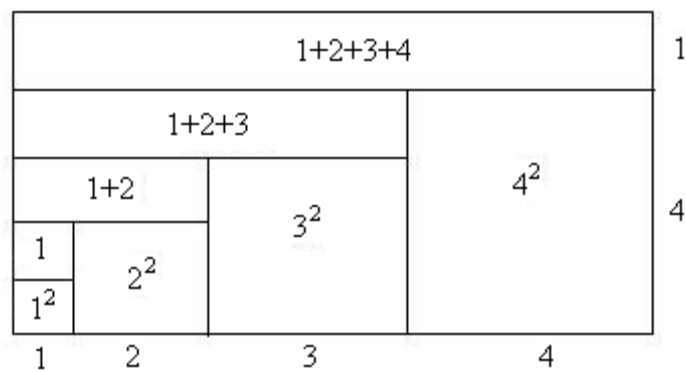

The area of this rectangle we can calculate in two different ways. The width is $\sum_{i=1}^{4}i$ and the height is $(4+1)$, i.e. the area is equal to $(4+1)\sum_{i=1}^{4}i$. On the other hand, the sum of the areas of the regions of the rectangle is equal to $\sum_{i=1}^{4}i^2 + \sum_{p=1}^{4}\left(\sum_{i=1}^{p}i\right)$.

We may generalize this equation replacing the number 4 with n:

$$(n+1)\sum_{i=1}^{n} i = \sum_{i=1}^{n} i^2 + \sum_{p=1}^{n}\left(\sum_{i=1}^{p} i\right),$$

where each side of the equation gives the area of the rectangle.

The diagram below illustrates the relationship

$$(4+1)\sum_{i=1}^{4} i^2 = \sum_{i=1}^{4} i^3 + \sum_{p=1}^{4}\left(\sum_{i=1}^{p} i^2\right).$$

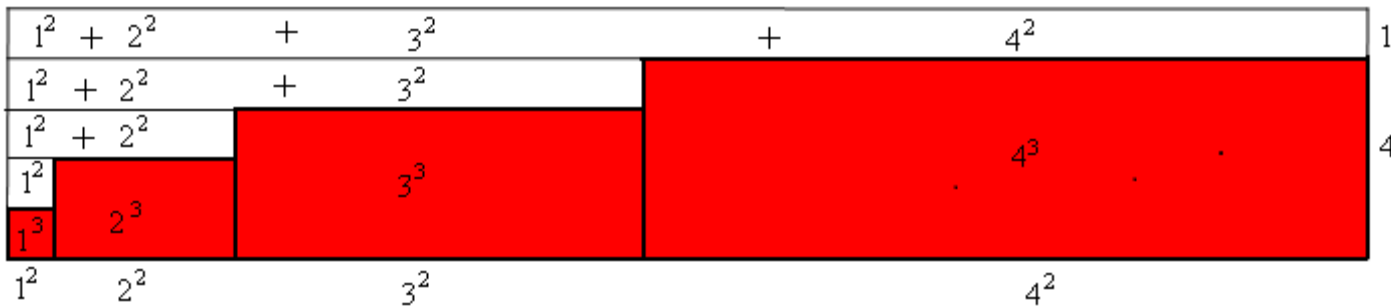

The diagram bellow illustrates the more general relationship.

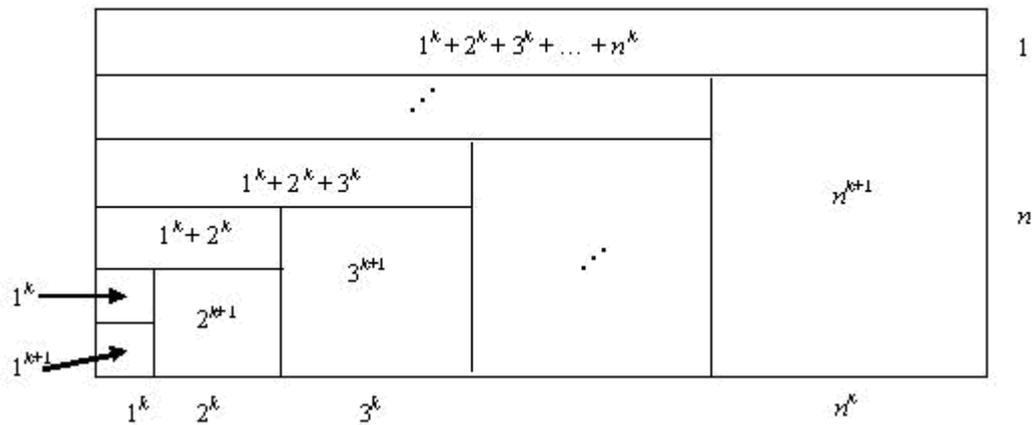

So we may conclude that

$$(n+1)\sum_{i=1}^{n} i^k = \sum_{i=1}^{n} i^{k+1} + \sum_{p=1}^{n}\left(\sum_{i=1}^{p} i^k\right)$$

Even more, we may put instead $i^k$ any function f(x).

For example if f(n)=n! we may derive the identity

$$(n+1)\sum_{i=1}^{n} i! = \sum_{i=1}^{n} i \cdot i! + \sum_{p=1}^{n}\left(\sum_{i=1}^{p} i!\right)$$

### 2.6. Abu Bakr ibn Muhammad ibn al Husayn al-Karaji (953 – 1029), Baghdad Iraq

Implicit proof by induction for arithmetic sequences was introduced by al-Karaji (c. 1000) and continued by **al-Samaw'al**, who used it for special cases of the binomial theorem and properties of Pascal's triangle.

He studied the Diophantine equation $x^3+y^3=z^2$ in area of the rational numbers.

Al-Karaji founded infinitely many solutions by parameterization

$$x = \frac{u^2}{1+v^3} \qquad y = \frac{u^2 v}{1+v^3} \qquad z = \frac{u^3}{1+v^3}.$$

He also shows the validity of identities

$$\sqrt[3]{54} - \sqrt[3]{2} = \sqrt[3]{16}$$

and

$$\sqrt[3]{54} + \sqrt[3]{2} = \sqrt[3]{128}.$$

Besides quadratic equations, he also considers the equations of higher degrees reduced to the square.

In this treatise, he gives more than 250 algebraic tasks for indeterminate equations. Al-Karaji is seen by many as the first person that completely free algebra from geometrical operations and to replace them with the arithmetical type of operations that are at the core of algebra today. He was first to define the monomials $x$, $x^2$, $x^3$, ... and $1/x$, $1/x^2$, $1/x^3$, ... and he gives rules for products of any two of these. He started a school of algebra, which flourished for several hundreds of years. This mathematician worked in Baghdad in Iraq and died in 1019. He proved the equality

$$1^3 + 2^3 + 3^3 + \cdots + 10^3 = (1 + 2 + \cdots + 10)^2.$$

In his argumentation, a contemporary reader will find an induction backward step.

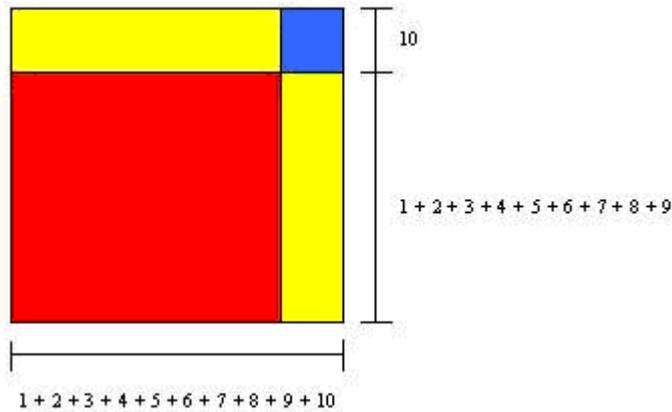

The area of the shown square is
$$S = (1+2+\cdots+10)^2.$$

The area of the two yellow rectangles and the blue square is
$$S_y + S_b = 2.10.(1+2+\cdots+9) + 10^2 = 2.10.\frac{10.9}{2} + 10^2 = 10^3$$

This gives us another way to write the area of the larger square namely as the sums of the areas of the red square and $S_y+S_b$.
$$S = S_r + S_y + S_b \Rightarrow (1+2+\cdots+9+10)^2 = (1+2+\cdots+9)^2 + 10^3.$$

There is the induction step: Repeating the preceding arguments for the red for the red square, we will have:
$$(1+2+\cdots+9)^2 = (1+2+\cdots+8)^2 + 9^3.$$

So
$$(1+2+\cdots+9+10)^2 = (1+2+\cdots+8)^2 + 9^3 + 10^3.$$

After repeating the procedure seven more times, we obtain the assertion. [8]

### 2.7. Ibn Tahir al-Baghdadi (980–1037) Baghdad

His name is connected with the Chinese Remainder Problem.
$$x \equiv 2(\bmod 3) \equiv 3(\bmod 5) \equiv 2(\bmod 7),$$

where x is an unknown that satisfies the requirements given in the remainder problem and needs to be determined. Formula of solution in the *Sun Zi Suanjing*:

$$70p+21r+15s-105n.$$

During the 11th century, mathematician Ibn Tahir al-Baghdadi discussed the Chinese Remainder Theorem in his treatise *AITakmilafi 'lim ai-Hisab*. The moduli that Ibn Tahir gave were the same as *Sun Zi Suanjing,* which were
$$m_1 = 3, \; m_2 = 5, \; m_3 = 7.$$

However, his problem was
$$x \equiv a \;(\bmod 3) \equiv b \;(\bmod 5) \equiv c \;(\bmod 7),$$

which was not entirely the same as *Sun Zi Suanjing.*

It was clear that Ibn Tahir had advanced further in his discussion of the remainder problem where arbitrary remainders, a, *b*, and *c* were given in his problem. It is interesting to note that Ibn Tahir was the first mathematician in antiquity to give an explanation regarding why the numbers 70, 21 and 15 were related to the moduli 3, 5 and 7 respectively. [3]

The popular Chinese Remainder Theorem found its way to Europe in a famous mathematical treatise by Italian mathematician Leonardo Fibonacci in 1202 entitled *Liber Abaci*. Even though the moduli given were the same, the remainder problem in *Liber Abaci* $x = 2 \pmod 3 = 3 \pmod 5 = 4 \pmod 7$, was different than the one in *Sun Zi Suanjing* in the sense that one of the remainders given was different. [7]

The contemporary explanation of the solution is based on the reciprocal number by modulus concept. If m and a are two mutually prime integers, then there exists only one in the set of residuals *x* such that

$$ax \equiv 1 \pmod{m}.$$

The solution of this modulo equation is denoted by $x^{-1}{}_{(m)}$.

If we have the system of equations by modulus:

$$\begin{vmatrix} x \equiv r_1 \pmod{m_1} \\ x \equiv r_2 \pmod{m_2} \\ \ldots \\ x \equiv r_n \pmod{m_n} \end{vmatrix}$$

then we denote by *M* the product $m_1 m_2 \ldots m_n$. By notation $M_i$ we mean the integer $M_i = M/m_i$. Then one solution of the system in the set of integers is

$$x_0 = r_1 M^{-1}_{(m_1)} M_1 + r_2 M^{-1}_{(m_2)} M_2 + r_3 M^{-1}_{(m_3)} M_3 + \cdots + r_n M^{-1}_{(m_n)} M_n.$$

In the integers, the system has infinitely many solutions given by formula

$$x = x_0 + tM.$$

The smallest one is given by

$$t = -\left\lfloor \frac{x_0}{M} \right\rfloor + 1.$$

**2.8. Omar Khayyam** (1048 –1131) Nishapur in North Eastern Iran, at a young age he moved to Samarkand

Omar Khayyam gave a complete classification of cubic equations with geometric solutions found by means of intersecting conic sections. Khayyam also wrote that he hoped to give a full description of the algebraic solution of cubic equations in a later work.

**2.9. Ibn Yahya al-Maghribi Al-Samawal** (1130-1180) Baghdad, Iraq and Maragha, Iran

Al-Samawal, nearly 200 years later, was an important member of al-Karaji's school. Al-Samawal was the first to give the new topic of algebra a precise description when he wrote that it was concerned: "... with operating on unknowns using all the arithmetical tools, in the same way as the arithmetician operates on the known."

**2.10. Sharaf al-Din al-Tusi** (c. 1135 – 1213) Tus, Iran, Aleppo, Mosul

al-Tusi although almost exactly the same age as al-Samawal, does not follow the general development that came through al-Karaji's school of algebra but rather follows Khayyam's application of algebra to geometry. He wrote a treatise on cubic equations, which: "... represents an essential contribution to another algebra which aimed to study curves by means of equations, thus inaugurating the beginning of algebraic geometry". Al-Tusi is recognized as the "father of trigonometry".

### 2.11. Hasan ibn al Qalasadi al-Andlaus (1412 – 1486)

Al-Qalaṣādī is taking the first steps toward the introduction of algebraic symbolism. In 1480 the Christian forces, raided and often pillaged the city, al-Qalasadi took refuge with his family in Tunisia, where he died in 1486. Baza was eventually besieged by the forces of Ferdinand and Isabella and its inhabitants sacked.